\documentclass[10pt]{article}

\usepackage{amsmath}
\usepackage{amsthm}
\usepackage{amssymb}
\usepackage{multicol}
\usepackage[active]{srcltx}
\usepackage{array}
\usepackage[dvipsnames,svgnames]{xcolor}
\usepackage{float}
\usepackage{pgfplots}
\usepgfplotslibrary{fillbetween}
\usepackage{pgfplots}
\pgfplotsset{compat=1.18}

\pgfplotsset{compat=1.18}
\usepackage{xcolor}
\definecolor{myred}{RGB}{180,0,0}

\usepackage{algorithm}
\usepackage{algpseudocode}

\usepackage{url}
\usepackage{tikz}
\usepackage{graphicx}
\usepackage{subcaption} 
\usepackage{pgfplots}
\usepackage{tabularx}
\usepackage{booktabs}
\usepackage[normalem]{ulem} 
\usepackage{comment}

\newtheorem{theorem}{Theorem}

\newtheorem{example}{Example}

\newtheorem{lemma}{Lemma}

\newtheorem{remark}{Remark}

\title{Subgradient and Stochastic Finite-Difference Methods for Constrained Strongly Quasiconvex Optimization}
\author{Rafik Khachatryan\thanks{Yerevan State University, 1 Alek Manukyan St, Yerevan 0025, Armenia. Email:khrafik@ysu.am, khachatryan.rafik@gmail.com}
\and
Felipe Lara\thanks{Instituto de Alta investigaci\'on (IAI), Universidad de Tarapac\'a, Arica, Chile. E-mail: felipelaraobreque@gmail.com; flarao@academicos.uta.cl. Web: felipelara.cl, ORCID-ID: 0000-0002-9965-0921}
}

\begin{document}

\maketitle
 
\begin{abstract}
\noindent
We propose a projected subgradient method for solving quasiconvex minimization problems. By employing the Clarke subdifferential, we prove that the sequence generated by the proposed algorithm converges to a solution of the minimization problem when the objective function is locally Lipschitz continuous and strongly quasiconvex in the sense of Polyak, under the same assumptions required for convex functions with the convex subdifferential. Moreover, we derive explicit convergence rate estimates for both the iterates and objective function values. In particular, compared with the corresponding results for nonsmooth, locally Lipschitz, strongly convex optimization, our analysis replaces strong convexity with strong quasiconvexity while preserving the order of the convergence rate, thereby extending the guarantees to a class of possibly nonconvex objective functions. In addition, we investigate a stochastic finite-difference scheme for minimizing a strongly quasiconvex function subject to inequality constraints defined by strongly quasiconvex functions, and establish a convergence rate of order $O(k^{-1/4})$, consistent with known complexity bounds for stochastic zeroth-order methods.
\end{abstract}

{\small}

\medskip


\medskip


\section{Introduction}
Modern convex and quasiconvex optimization rests on classical existence and convergence results, from Polyak's foundational work on minimizing sequences under constraints \cite{polyak1966} to Rockafellar's monotone operator theory and proximal point algorithm \cite{rockafellar1976}. Clarke's theory of generalized differentiation for locally Lipschitz functions \cite{clarke1975} later provided the tools needed to analyze subgradient and projection-type algorithms in nonsmooth settings. Vial \cite{vial1982} and Jovanovic \cite{jovanovic1996} clarified the structural differences between convex, strongly convex, and quasiconvex functions, while Nesterov's work on nonsmooth convex and quasiconvex minimization \cite{nesterov1984, nesterov1989, nesterov2010} introduced algorithmic principles later extended by Nesterov and Spokoiny \cite{nesterov2017} to derivative-free, gradient-free (two-point finite-difference) methods with optimal complexity guarantees via Gaussian smoothing.

More recently, attention has turned to strongly quasiconvex functions, a class bridging convexity and quasiconvexity while retaining useful stability and convergence properties. A series of works \cite{lara2024, hadjisavvas2026a, hadjisavvas2026b, KL-2, lara2022, lara2025a} established existence results, characterizations, and algorithmic schemes (proximal point and subgradient projection methods) for quasiconvex and strongly quasiconvex problems, including convergence via strong subdifferentials \cite{KL-2} under assumptions analogous to the convex case, and a quasi-linear convergence rate for the iterates. Related dynamical-system interpretations and gradient-based methods appear in \cite{hadjisavvas2026b, lara2025a}; see the survey \cite{grad2025} for a fuller account.

Stochastic and gradient-free methods for quasiconvex optimization have also gained recent attention. Khachatryan \cite{Khachatryan2026Vestnik} studied stochastic finite-difference methods for quasiconvex black-box optimization with a maximum-type objective $f(x) = \max_{i\in I} f_i(x)$, where each $f_i$ is quasiconvex with Lipschitz-continuous gradient on a convex compact set, obtaining an $O(1/\sqrt{k})$ convergence rate. However, the feasible set there is a plain set constraint $x \in M$; no inequality constraints $f_i(x)\le 0$ with quasiconvex $f_i$ are considered, so constrained quasiconvex optimization and structure-preserving penalty functions remain unaddressed.

Zeroth-order methods for functionally constrained optimization have also been developed, notably by Nguyen and Balasubramanian \cite{NguyenBalasubramanian2022} and by Li et al. \cite{LiChenLiuLuXu2022}, who proposed a zeroth-order inexact augmented Lagrangian method for composite problems with functional constraints. These works motivate studying stochastic zeroth-order methods in the quasiconvex setting considered here.

Finally, Gasnikov et al.~\cite{gasnikov2024} recently studied smooth quasiconvex optimization in the comparison-oracle model, obtaining complexity guarantees from function-value comparisons alone, situating quasiconvex optimization within the broader family of oracle-based and zeroth-order methods.

\paragraph{Our contributions.}
This paper broadens the convergence theory of stochastic finite-difference methods from convex to quasiconvex optimization. This is not merely a formal weakening of assumptions: quasiconvexity includes nonconvex, nonsmooth functions lying outside the classical convex framework, and since the convex subgradient inequality is no longer available, the analysis instead exploits geometric relations between generalized gradients, distances to optimal solutions, and the structure of quasiconvex and strongly quasiconvex functions. On this basis, we develop a convergence theory for stochastic finite-difference methods on locally Lipschitz quasiconvex functions, covering problems outside the standard convex setting, and we establish convergence results and rates for the strongly quasiconvex case; since every strongly convex function is strongly quasiconvex but not conversely, this framework strictly contains the classical strongly convex setting. We further extend the analysis to constrained quasiconvex optimization through a penalty-function method for problems with strongly quasiconvex functional constraints, proving convergence under suitable regularity and nondegeneracy assumptions even when the constraints are nonconvex. Throughout, local Lipschitz continuity and these nondegeneracy conditions stand in for the missing convex structure, and the resulting rates give quantitative guarantees for a substantially broader, genuinely nonconvex class of problems than the classical convex theory covers.

\medskip

The paper is organized as follows. Section \ref{sec:02} recalls preliminaries on generalized convexity and nonsmooth analysis. Section \ref{sec:03} presents a subgradient method for locally Lipschitz strongly quasiconvex functions based on the Clarke subdifferential, with a complete convergence analysis for both iterates and function values. Section \ref{sec:04} develops a stochastic finite-difference method for a mathematical programming problem whose constraint functions are locally Lipschitz and strongly quasiconvex and whose objective is locally Lipschitz and quasiconvex, establishing a convergence rate of $O(k^{-1/4})$ under mild assumptions. Section \ref{sec:05} develops a penalty function method for the same problem under quasiconvexity alone, giving a sufficient condition for convergence of the generated sequence to the optimal solution. Numerical illustrations support our theoretical findings.

\section{Preliminaries and Basic Definitions} \label{sec:02}

The inner product of $\mathbb{R}^{n}$ and the Euclidean norm are denoted by $\langle \cdot,\cdot \rangle$ and $\lVert \cdot \rVert$, respectively. The closed  ball with center at $x_{0} \in \mathbb{R}^{n}$ and radius $\varepsilon > 0$ is denoted by $B_{\varepsilon} (x_{0})$. 



Given any extended-valued function $h: \mathbb{R}^{n} \rightarrow \overline{\mathbb{R}}$, the effective domain of $h$ is defined by ${\rm dom}\,h := \{x \in \mathbb{R}^{n}: h(x) < + \infty \}$. 

A function $h$ with convex domain is said to be:
\begin{itemize}
 \item[$(a)$] (strongly) convex on ${\rm dom}\,h$, if there exists $\mu \geq 0$ such that, for all $x, y \in {\rm dom}\,h$ and all $\lambda \in [0, 1]$, we have
 \begin{equation}\label{strong:convex}
  h(\lambda y + (1-\lambda)x) \leq \lambda h(y) + (1-\lambda) h(x) - \lambda (1 - \lambda) \frac{\mu}{2} \lVert x - y \rVert^{2},
 \end{equation}

 \item[$(b)$] (strongly) quasiconvex on ${\rm dom}\,h$, if there exists $\mu \geq 0$ such that, for all $x, y \in {\rm dom}\,h$ and all $\lambda \in [0, 1]$, we have
 \begin{equation}\label{strong:quasiconvex}
  h(\lambda y + (1-\lambda)x) \leq \max \{h(y), h(x)\} - \lambda(1 - \lambda) \frac{\mu}{2} \lVert x - y \rVert^{2}.
 \end{equation}
\end{itemize} 
 \noindent A function is strongly convex (resp. quasiconvex) when $\mu > 0$, and convex (resp. quasiconvex) when $\mu = 0$. Hence, every (strongly) convex function is (strongly) quasiconvex, while the reverse statements does not hold in general (see \cite{CMA,HKS,lara2022}).

A locally Lipschitz function $f$ is called regular in the sense of Clarke on $M$ if, for every $x\in M$ and every $v\in\mathbb{R}^n$, the ordinary
directional derivative exists and satisfies
\[
 f'(x;v)=f^\circ(x;v).
\]


\begin{remark}\label{rem:one}
\begin{itemize}
 \item[$(i)$] Let us construct an example of a strongly quasiconvex function that is locally Lipschitz, not differentiable at all points, and not convex. Indeed, let $f: \mathbb{R} \rightarrow \mathbb{R}$ be the function given by
\[
f(x)=\max\{f_1(x), f_2(x)\},
\]
where $f_1(x) = \max\{0, x\} + \lvert x \rvert^2$ and $f_2 (x) = x^2 + 2 \sin^{2} (x)$.

Here, $f_1$ is strongly convex, thus strongly quasiconvex with modulus $\mu_{1} = 2$. The function \(f_2\) is also strongly quasiconvex with modulus $\mu_{2} = \frac{1}{2}$ (see \cite[Example 30]{lara2025a}).

Since the maximum of strongly quasiconvex functions is again strongly quasiconvex \cite[Corollary 2]{VNC-2}, it follows that \(f\) is strongly quasiconvex on the whole domain with modulus $\mu := \min\{\mu_{1}, \mu_{2}\} = \frac{1}{2}$.

On the other hand, note that $f$ is locally Lipschitz as the maximum of two locally Lipschitz functions, it is not differentiable at points where $f_1(x)=f_2(x)$, and it is not convex due to $f_2$.

 \item[$(ii)$] Let $h: \mathbb{R}^{n} \rightarrow \mathbb{R}$ be a strongly convex (possible nonsmooth) and positive function. Let us show that the function
\[
 f(x) = h(x) - \frac{1}{h(x) + a}, \quad a > 0,
\]
is strongly quasiconvex. Indeed, let us define the functions:
$$h_1(x)=h(x), \qquad h_2 (x) = - \frac{1}{h(x)+a}.$$
Note that $h_2$ is quasiconvex since its sublevel sets are convex. Now, take $x$, $y$, and define $z := (1 - \lambda)x + \lambda y$ with $\lambda \in [0, 1]$. Suppose without loss of generality that $h_1(x) \leq h_1(y)$. This implies $h_2(x) \leq h_2(y)$. By quasiconvexity of $h_2$, we get $h_2(z) \leq h_2(y)$.
Thus,
\[
\begin{aligned}
 h(z) &= h_1(z) + h_2(z) \\
 & \leq \max\{h_1(x), h_1(y)\} + h_2(y) - \frac{\mu}{2} \lambda(1 - \lambda)\|y - x\|^2 \\
 & \leq h(y) - \frac{\mu}{2}\lambda(1 - \lambda) \|y - x\|^2.
\end{aligned}
\]
Therefore, $h$  is strongly quasiconvex with modulus $\mu>0$.
\end{itemize}
\end{remark}

\section{A Subgradient Method}\label{sec:03}

This section develops the auxiliary results needed to establish the convergence properties of the proposed method. We first introduce the technical quantity that will be used throughout the analysis and explain its role in the subsequent estimates. We then present several auxiliary lemmas that provide the key inequalities required in the proof of the main convergence result. 
For clarity, we indicate after each result how it is used in the subsequent analysis.

Given a convex set $M \subset \mathbb{R}^{n}$, for arbitrary $x, y \in M \subset \mathbb{R}^{n}$ and $g(x) \in \partial f(x)$ with $g(x) \neq 0$, we introduce the auxiliary quantity
\begin{equation}\label{eq:02}
 v_{f}(x,y) := \left\langle \frac{g(x)}{\|g(x)\|},x - y \right \rangle. 
\end{equation}

A fundamental inequality used in our convergence analysis is provided below. It relates the subgradient of a strongly quasiconvex function at an arbitrary point $x\in M$ to the distance from that $x$ to the unique minimizer $x^*$. 


\begin{lemma} \label{lemma:03}
 Let \(M\subset\mathbb{R}^{n}\) be a closed convex set, and let
\(f:\mathbb{R}^{n}\to\mathbb{R}\) be a locally Lipschitz continuous,
regular, and strongly quasiconvex function with modulus \(\mu>0\) on \(M\). Then, for all $x \in M$, we have
\[
 \frac{\mu}{2} \| x - x^{*}\|^{2} \leq \langle v,x - x^{*} \rangle, ~ \forall ~ v \in \partial f(x),
\]
where \(x^{*}\) is the unique minimizer of \(f\) on \(M\). 
\end{lemma}

\begin{proof}
Since $x^*$ is a global minimizer of $f$ on $M$, we have $f(x^*)\leq f(x)$ for all $x\in M$. Using this and the strong quasiconvexity of $f$, we have for every $t\in[0,1]$ that
\[
 f(x+ t(x^*-x)) \leq f(x) -\frac{\mu}{2} t(1-t) \|x-x^*\|^2.
\]
Dividing by $t>0$, we obtain
\[
 \frac{ f(x+t (x^*-x))-f(x)}{t} \leq - \frac{\mu}{2} (1-t) \|x-x^*\|^2.
\]
Letting $t\downarrow0$ and using the regularity of $f$, which implies that the Clarke directional derivative coincides with the ordinary directional derivative, we obtain
\[
 f^\circ(x;x^*-x)=f'(x;x^*-x) \leq -\frac{\mu}{2}\|x-x^*\|^2.
\]

For a locally Lipschitz function, the Clarke directional derivative satisfies $f^\circ (x; d) = \max_{v \in\partial f(x)} \langle v, d \rangle$. Taking $d=x^*-x$, we have
\[
 \max_{v \in \partial f(x)} \langle v, x^* - x \rangle \leq -\frac{\mu}{2} \|x-x^*\|^2 ~ \Longrightarrow ~ \langle v,x-x^*\rangle \geq \frac{\mu}{2} \|x-x^* \|^2, ~ \forall ~v \in \partial f(x).
\]
This proves the lemma.
\end{proof}

\begin{remark}
The regularity assumption cannot be omitted. Indeed, let
$M=[0,1]$ and $x^*=1$, and define
\[
f(x)=
\begin{cases}
(x-1)^2+102x, & x<0,\\
(x-1)^2, & x\geq 0.
\end{cases}
\]
Then $f$ is locally Lipschitz on $\mathbb{R}$ and is $2$-strongly quasiconvex on $M$. However,
\[
\partial f(0)=[-2,100].
\]
 Taking $v = 100 \in \partial f(0)$, we obtain $\frac{\mu}{2} \lvert 0-1 \rvert^2 = 1 > 100(0-1) = - 100$. Thus, the claimed inequality may fail without the regularity assumption. Indeed, $f$ is not Clarke regular at $x=0$, since $f'(0;1) = -2$ while $f^\circ (0; 1) = 100$. Hence, $f'(0; 1) \neq f^\circ (0; 1)$, which shows that $f$ is not regular at $0$.
\end{remark}

The following lemma may be viewed as a quasiconvex analogue of the
first-order inequality in the convex case. It provides a subgradient-based
estimate for the functional-value difference and will be used later in the
analysis of the penalty function method.

\begin{lemma}\label{lemma:04} 
 Let \(M\subset \mathbb{R}^{n}\) be a convex compact set, and $f:\mathbb{R}^{n} \to \mathbb{R}$ be a strongly quasiconvex and locally Lipschitz function on \(\mathbb{R}^{n}\). Then for every point \(x, y \in M\) with \(f(x) \geq f(y)\), the following inequality holds:
\begin{equation}\label{eq:03}
 f(x) - f(y) \leq L_{\Omega} v_{f} (x, y), ~ \forall ~ g \in \partial f(x)
\end{equation}
where \(L_{\Omega}\) is the Lipschitz constant of the function \(f\) on the set $\Omega=M + {\rm diam}(M) B_{1}(0)$. 
\end{lemma}

\begin{proof}
 Let $z = y + \| g\|^{-2} \langle g, x - y\rangle g$. Clearly, \( \langle g, z - x \rangle = 0\). Hence, $f(z) \geq f(x)$. Indeed, suppose for the contrary that $f(z)< f(x)$. Hence,
\[
 f(x + \lambda (z - x)) \leq \max \{f(z), f(x)\} -\frac{\mu}{2} \lambda (1 - \lambda)\| x - z\|^{2}.
\]
It follows that,
\[
\begin{aligned}
 & f^{\circ} (x,z - x) = \limsup_{u \to x,\lambda \to 0}\frac{f(u + \lambda (z - x)) - f(u)}{\lambda} \\
 & = \limsup_{u\to x,\lambda \to 0}\frac{f(u + \lambda (z - u)) - f(u)}{\lambda} + L_{\Omega} \limsup_{u\to x}\| u - x\| - \frac{\mu}{2} (1 - \lambda)\| u - z\|^{2} \\
 & \leq - \frac{\mu}{2} \| z - x\|^{2}.
\end{aligned}
\]
From this it directly follows that for any \(g\in \partial f(x)\), we have $\langle g,z - x\rangle < 0$, a contradiction. Therefore, $f(z) \geq f(x)$.

Since $f$ is locally Lipschitz, it satisfies the global Lipschitz condition on the convex compact set $\Omega.$ Since the points \(x,y,z\) lie in this domain, it follows that
\[
 f(x) - f(y) \leq f(z) - f(y) \leq L_{\Omega}\| y - z\| = L_{\Omega} \frac{\langle g, x - y \rangle}{\|g\|} = L_{\Omega} v_{f}(x,y),
\]
and the result follows.
\end{proof}

The analysis in this section includes both strongly quasiconvex and quasiconvex minimization over convex compact sets. Note that for the quasiconvex case, the convergence analysis of the functional values requires a different, but closely related, subgradient property. The following lemma provides the corresponding quasiconvex counterpart of the inequality established above for strongly quasiconvex functions.

\begin{lemma}\label{coro:06}
{\rm (\cite[Corollary 1.1]{Khachatryan2026Vestnik})}
Let \(f:\mathbb{R}^{n}\to\mathbb{R}\) be a locally Lipschitz and quasiconvex
function on the set
\[
\Omega = M + {\rm diam}(M) B_{1}(0),
\]
and regular on \(M\). Then, for every \(x,y\in M\) such that
\(f(y)\leq f(x)\), and for every \(g(x)\in\partial f(x)\) with
\(g(x)\neq 0\), the following inequality holds:
\begin{equation}\label{eq:05}
f(x)-f(y)
\leq
L_{\Omega}
\frac{\langle g(x),x-y\rangle}{\|g(x)\|}.
\end{equation}
\end{lemma}

\begin{remark}
 The regularity assumption in Lemma \ref{coro:06} is essential for relation \eqref{eq:05}. Indeed, let $M = [0,1]$ and $f(x) = \min\{x, 0\}$. Then, $f$ is quasiconvex and locally Lipschitz, but it is not regular at $x=0$ since
\[
 f^\circ (0;d) = \max\{0,d\} \neq \min\{0, d\} = f'(0; d) \quad (d>0).
\]
 Moreover, $\partial f(0) = [0, 1]$. Taking $x=0$, $x^* = 1$, and $g(0) = 1 \in \partial f(0)$, the right-hand side of \eqref{eq:05} equals $-1$, whereas $f(0) - f(1) = 0$. Hence \eqref{eq:05} fails, showing that the regularity assumption cannot be omitted.
\end{remark}

Let us consider the following method for minimizing strongly quasiconvex functions on a closed convex set \(M\),
\[ 
 x^{k + 1} = \Pi_{M}(x^{k} - \rho_{k}v^{k}), ~ \forall ~ v^{k} \in \partial f(x^{k}),
\]
where \(\Pi_{M}\) denotes the projection onto the set \(M\) and \(\rho_{k} > 0\) is the step size.

\medskip

The following theorem provides the convergence rate of the iterates for the subgradient method applied to a strongly quasiconvex and locally Lipschitz function. It estimates how quickly the sequence of iterates $\{x^{k}\}_{k}$ approaches the unique minimizer \(x^{*}\) of \(f\) on \(M\).

\begin{theorem}\label{theo:01} 
 Let \(M\) be a compact convex set, and suppose that the function \(f\) satisfies all the assumptions of Lemma \ref{lemma:03}. Take \(\rho_{k} = \frac{1}{k}\) for all $k$. Then the following estimates hold:
$$
\|x^k-x^*\|^2 =
\begin{cases}
 \mathcal{O} \left(\dfrac{1}{k} \right), & \mu > 1, \\[1mm]
 \mathcal{O} \left( \dfrac{\ln k}{k} \right), & \mu = 1, \\[1mm]
 \mathcal{O} \left( \dfrac{1}{k^\mu} \right), & 0 < \mu < 1.
\end{cases}
$$
\end{theorem}

\begin{proof}
Observe that
\[ 
 \| x^{k + 1} - x^{*}\|^{2} \leq \| x^{k} - x^{*} \|^{2} + 2 \rho_{k} \langle v^{k},x^{*} - x^{k} \rangle + \rho_{k}^{2}\| v^{k}\|^{2}.
\]
Since \(x^{*}\) is a minimizer, Lemma \ref{lemma:03} implies
\[
\langle v^{k},x^{*} - x^{k}\rangle \leq 0, ~ \forall ~ v^{k} \in \partial f(x^{k}).
\]
Since the sequence $\{x^{k}\}_{k}$ lies in the compact set \(M\), it is bounded. Hence, the sequence \(\{v^{k}\}_{k}\) is also bounded, i.e., there exists a constant \(C > 0\) such that $\| v^{k}\| \leq C$ for all $k \in \mathbb{N}$. Using Lemma \ref{lemma:03}, we obtain
\begin{equation}\label{eq:07}
 \| x^{k + 1} - x^{*}\|^{2} \leq \left(1 - \frac{\mu}{k} \right)\| x^{k} - x^{*}\|^{2} + \frac{C^{2}}{k^{2}}. \tag{7}
\end{equation}

Denote \(a_{k} = \| x^{k} - x^{*}\|^{2}\). Then the recurrence relation is
\[
 a_{k + 1} \leq \left(1 -  \frac{\mu}{k} \right)  a_{k} + \frac{C^{2}}{k^{2}}.
\]

Iterating this inequality, we get
\[
 a_{k + 1} \leq a_{k_{0}} \prod_{i = k_{0}}^{k} \left(1 - \frac{\mu}{i}\right) + \sum_{j = k_{0}}^{k} \frac{C^{2}}{j^{2}} \prod_{i = j + 1}^{k} \left(1 - \frac{\mu}{i} \right) + \frac{C^{2}}{k^{2}}.
\]

Since \(\ln (1 - \mu /i) \leq - \mu /i\), we have
\[ 
 \prod_{i = j + 1}^{k}\left(1 - \frac{\mu}{i} \right) = \exp \Big(\sum_{i = j + 1}^{k}\ln (1 - \mu /i)\Big) \leq \exp \Big(-\mu \sum_{i = j + 1}^{k} \frac{1}{i} \Big).
\]

Also,
\[
 \sum_{i = j + 1}^{k}\frac{1}{i} \geq \int_{j + 1}^{k + 1} \frac{d x}{x} = \ln (k + 1) - \ln (j + 1),
\]
thus,
\[
 \prod_{i = j + 1}^{k} \left(1 - \frac{\mu}{i} \right) \leq \left( \frac{j + 1}{k + 1} \right)^{\mu} ~~~ {\rm and} ~~~
 \prod_{i = k_{0}}^{k} \left(1 - \frac{\mu}{i} \right) \leq \left(\frac{k_{0}}{k + 1} \right)^{\mu}.
\]

Therefore, the sum term can be estimated as
\[
 S_{k} = \sum_{j = k_{0}}^{k} \frac{1}{j^{2}} \prod_{i = j + 1}^{k} \left(1 - \frac{\mu}{i} \right) \leq \sum_{j = k_{0}}^{k} \frac{1}{j^{2}} \left( \frac{j + 1}{k + 1} \right)^{\mu} \leq \frac{C^{2}2^{\mu}}{(k + 1)^{\mu}} \sum_{j = k_{0}}^{k}j^{\mu -2}.
\]

Consider different cases of \(\mu\). If \(\mu < 1\), then
\[
 \sum_{j = k_{0}}^{k}j^{\mu -2} \leq \int_{k_{0}}^{k} x^{\mu -2}d x = \frac{k^{\mu -1} - k_{0}^{\mu -1}}{\mu - 1} \leq \frac{k_{0}^{\mu -1}}{1 - \mu}.
\]

Thus,
\[
 S_{k} \leq \frac{2^{\mu}k_{0}^{\mu -1}}{1 - \mu} \frac{1}{(k + 1)^{\mu}}.
\]
Therefore, in this case we obtain the asymptotic estimate $a_{k} = \mathcal{O} \left(\frac{1}{k^{\mu}} \right)$.

Similarly, we obtain $a_k = \mathcal{O} \left(\frac{\ln k}{k}\right)$ for $\mu = 1$, and  $a_k = \mathcal{O} \left(\frac{1}{k}\right)$ for $\mu > 1$, which completes the proof.
\end{proof}

\begin{remark}[Practical use of the strong quasiconvexity modulus] \label{rem:mu}
The modulus $\mu$ is a theoretical regularity parameter and need not be
known exactly in order to implement the proposed algorithm. If only a
conservative lower estimate $\widehat{\mu}$ is available, with
\[
0<\widehat{\mu}\leq\mu,
\]
then the convergence analysis remains valid with $\widehat{\mu}$ in
place of $\mu$. The resulting estimates may be more conservative than
those obtained with the exact modulus, but the algorithm itself is
unchanged. In particular, the exact value of $\mu$ is not required for
the implementation of the projected subgradient iteration.
\end{remark}

\begin{remark}[Polyak adaptive stepsize]
\label{rem:Polyak}
Throughout this remark we assume, as in Theorem~\ref{theo:01}, that $f$
is strongly quasiconvex with modulus $\mu>0$, and, in addition, that
$\gamma:=\min_{x\in M}\operatorname{dist}(0,\partial f(x))>0$ (the
nondegeneracy condition of Theorem~\ref{theo:02}).
The deterministic stepsize
\[
\rho_k=\frac{1}{k}
\]
used in Theorem~\ref{theo:01} is convenient for the theoretical
analysis and leads to explicit convergence rates. In practical
implementations, however, an adaptive Polyak-type stepsize can also be
considered. More precisely, let $\Delta_k:=f(x^k)-f(x^*)$ and consider
$\rho_k = \lambda\frac{\Delta_k}{\|v^k\|^2}$ with $v^k\in\partial
f(x^k)$, where $\lambda>0$ is a relaxation parameter. The projected
subgradient iteration then takes the form
\[
x^{k+1}
=
\Pi_M\left(
x^k-
\lambda
\frac{f(x^k)-f(x^*)}{\|v^k\|^2}v^k
\right).
\]
The main difference from the prescribed stepsize $\rho_k=1/k$ is that
the Polyak stepsize adapts automatically to the current objective
error: larger steps are taken when the iterate is far from the
solution, and the stepsize decreases automatically as $f(x^k)-f(x^*)$
approaches zero.

Indeed, using the projection inequality and the quasiconvexity
estimate
\[
\langle v^k,x^k-x^*\rangle \geq \frac{\gamma}{L_\Omega}
\bigl(f(x^k)-f(x^*)\bigr),
\]
we obtain
\[
\begin{aligned}
\|x^{k+1}-x^*\|^2
& \leq \|x^k-x^*\|^2 -2\rho_k\langle v^k,x^k-x^*\rangle+\rho_k^2\|v^k\|^2\\
& \leq \|x^k-x^*\|^2 - \lambda\left(\frac{2\gamma}{L_\Omega}-\lambda\right)
\frac{(f(x^k)-f(x^*))^2}{\|v^k\|^2}.
\end{aligned}
\]
Consequently, the stepsize is descent-producing whenever
$0<\lambda<\frac{2\gamma}{L_\Omega}$.

If, in addition, we use the quadratic growth property of strongly
quasiconvex functions (see \cite[Corollary 9]{hadjisavvas2026a}),
\[
f(x)-f(x^*)\geq\frac{\mu}{4}\|x-x^*\|^2,
\]
together with the uniform bound $\|v^k\|\leq C$ on the subgradients,
then
\[
\|x^{k+1}-x^*\|^2 \leq \|x^k-x^*\|^2 - c\|x^k-x^*\|^4
\]
for some $c>0$, and hence $\|x^k-x^*\|^2=\mathcal O(1/k)$.

Thus, the Polyak stepsize provides an adaptive alternative to the
deterministic choice $\rho_k=1/k$: it uses the current
objective-function error rather than the iteration counter to
determine the stepsize. As in the classical Polyak scheme, the exact
optimal value $f(x^*)$ is required; if it is unavailable, an estimate
or a target value can be used instead.
\end{remark}

We next establish a functional-value convergence rate for the projected
method when $f$ is quasiconvex, without assuming strong quasiconvexity.

\begin{theorem}\label{theo:02}
 Let \(M\) be a compact convex set, and suppose that the function $f$ satisfies all the assumptions of Lemma \ref{coro:06}. If, in addition,
\begin{equation} \label{gamma:01}
 \gamma \equiv \min_{x\in M} {\rm dist} \bigl(0, \partial f(x) \bigr)>0,
\end{equation}
then, for the projected method with $\rho_k = \frac{1}{\sqrt{k}}$ for all $k$, the following estimate holds:
\[
 \min_{1 \leq k \leq T} \bigl( f(x^k) - f(x^*) \bigr) = \mathcal{O} \left(\frac{1}{\sqrt{T}} \right).
\]
\end{theorem}

\begin{proof}
 Since $x^k, x^* \in M$ and $M$ is compact, we have, for every $k$, $\|x^k-x^*\| \leq {\rm diam} (M)$, and hence $\|x^k-x^*\|^2 \leq {\rm diam}^2 (M)$. Furthermore, from the projection step  $x^{k+1} = \Pi_M (x^k-\rho_k v^k)$, the nonexpansiveness of the projection gives
\[
 \|x^{k+1}-x^*\|^2 \leq \|x^k-x^*\|^2 -2 \rho_k \langle v^k, x^k-x^*\rangle + \rho_k^2 \|v^k\|^2.
\]
Using inequality \eqref{eq:05} and the assumption $\|v^k\|\geq\gamma$, we obtain
\begin{equation}\label{eq:08_new}
 f(x^k)-f(x^*) \leq \frac{L_{\Omega}}{2 \gamma} \left( \frac{\|x^k-x^*\|^2 - \|x^{k+1}-x^*\|^2}{\rho_k} + C^2\rho_k \right).
\end{equation}

Summing \eqref{eq:08_new} for $k=1,\ldots,T$ gives
\begin{equation}\label{eq:09_new}
\begin{aligned}
 \sum_{k=1}^{T} (f(x^k)-f(x^*)) \leq \frac{ L_{\Omega} }{2\gamma} \Bigg[ & \sum_{k=1}^{T} \|x^k - x^*\|^2 \left( \frac{1}{\rho_k}-\frac{1}{\rho_{k-1}} \right) + C^2\sum_{k=1}^{T}\rho_k \Bigg],
\end{aligned}
\end{equation}
where $1/\rho_0=0$. Set $\rho_k=\frac{1}{\sqrt{k}}$. Then $\frac{1}{\rho_k}=\sqrt{k}$ and $\frac{1}{\rho_k}-\frac{1}{\rho_{k-1}} = \sqrt{k}-\sqrt{k-1}$. Thus,
\[
\begin{aligned}
 \sum_{k=1}^{T} \|x^k-x^*\|^2 \left( \frac{1}{\rho_k} - \frac{1}{\rho_{k-1}} \right)
 & \leq {\rm diam}^2 (M) \sum_{k=1}^{T} \left( \sqrt{k} -\sqrt{k-1}\right) \\
&=
\operatorname{diam}^2(M)\sqrt{T}.
\end{aligned}
\]

On the other hand, $\sum_{k=1}^{T}\rho_k = \sum_{k=1}^{T}\frac1{\sqrt{k}} \leq 2\sqrt{T}$. Substituting these estimates into \eqref{eq:09_new}, we obtain
\[
 \sum_{k=1}^{T} (f(x^k)-f(x^*)) \leq \frac{L_\Omega}{2 \gamma} ( {\rm diam}^2 (M) + 2C^2) \sqrt{T}.
\]
Since $T\min_{1\leq k\leq T} ( f(x^k)-f(x^*) ) \leq
\sum_{k=1}^{T} (f(x^k)-f(x^*))$, we conclude that
\begin{align*}
 & \hspace{0.8cm} \min_{1 \leq k \leq T} (f(x^k) - f(x^*)) \leq \frac{L_{\Omega}}{2 \gamma} \frac{ {\rm diam}^2 (M) + 2 C^2}{\sqrt{T}} \\
 & \Longrightarrow ~ \min_{1\leq k\leq T} (f(x^k) -f(x^*)) = \mathcal{O} \left( \frac1{ \sqrt{T}} \right),
\end{align*}
which completes the proof.
\end{proof}

\begin{remark}
 The assumption $\gamma>0$ cannot be removed in general. Indeed, consider $M=[-1,1]$ and $f(x)=x^3$. Then $f$ is quasiconvex but not convex, with unique minimizer $x^*=-1$, whereas
\[
 \gamma = \min_{x \in M} {\rm dist} (0, \partial f(x)) = \min_{x \in [-1, 1]} 3x^2 = 0.
\]
Starting from $x^0=0$, the projected subgradient method with $\rho_k=(k+1)^{-1/2}$ remains at $x^k =0$, since $0\in\partial f(0)$. Hence $f(x^k)-f(x^*)=1$ for all $k$, and no convergence to the optimal value occurs. This shows that $\gamma>0$ is essential for the convergence result of Theorem \ref{theo:02}.
\end{remark}

\begin{remark}[Relation to convex and strongly convex optimization]
The results of this section extend the standard nonsmooth subgradient
framework for convex optimization to the broader class of locally
Lipschitz quasiconvex functions. As shown in \cite{Karimi2016}, when
$f$ is strongly convex with modulus $\mu$ and has an $L$-Lipschitz
gradient, gradient descent with step size $1/L$ converges linearly,
\[
f(x^k)-f(x^*)\leq\left(1-\frac{\mu}{L}\right)^k\bigl(f(x^0)-f(x^*)\bigr).
\]
Theorems \ref{theo:01} and \ref{theo:02} extend this circle of ideas to a substantially broader class of functions (quasiconvex and locally Lipschitz, for which the gradient may not exist) at the cost of sublinear rather than linear rates. This weakening is unavoidable: even for merely convex, nonsmooth objectives, subgradient methods cannot in general achieve better than sublinear rates, so the results obtained here are consistent with, rather than weaker than, what smoothness gives up.

The extension is nontrivial also at the level of the underlying inequalities, since the fundamental convex subgradient inequality
\[
 f(x)-f(y) \leq \langle g(x), x-y \rangle,\qquad g(x) \in \partial f(x),
\]
is generally not valid for quasiconvex functions. Instead, under the regularity assumption, we employ the first-order characterization
\[
 \langle g(x),y-x\rangle\leq0,\qquad f(y)\leq f(x),
\]
which, together with local Lipschitz continuity, yields the normalized subgradient estimate
\[
 f(x)-f(x^*) \leq L_\Omega \left\langle \frac{g(x)}{\|g(x)\|}, x-x^* \right\rangle, \qquad g(x) \in \partial f(x),\ g(x)\neq0.
\]
Thus, the quasiconvex setting requires a different first-order inequality and, in general, an additional nondegeneracy condition $\gamma>0$.

For strongly quasiconvex functions, the stronger estimate
\[
\frac{\mu}{2}\|x-x^*\|^2\leq\langle g,x-x^*\rangle,\qquad g\in\partial f(x),
\]
provides direct control of the distance to the minimizer without
requiring $\gamma>0$: under the assumptions of Theorem~\ref{theo:01},
this condition is not needed, whereas it is essential in
Theorem~\ref{theo:02} to control the normalized subgradient.

Consequently, the passage from convexity to quasiconvexity changes not
only the class of admissible objective functions but also the
analytical tools and assumptions required for convergence: the
standard convex subgradient inequality is replaced by a normalized
quasiconvex estimate, and the regularity assumption is used to connect
the Clarke subdifferential with directional derivatives. This confirms
that the obtained extension is substantive, and not merely a formal
replacement of convexity by quasiconvexity.
\end{remark}

\section{A Stochastic Finite-Difference Method}
\label{sec:04}

In this section, we develop a projected stochastic finite-difference method for constrained quasiconvex optimization and establish explicit convergence estimates for both the objective error and constraint violation. The method employs two-point randomized finite-difference estimators in place of exact subgradient evaluations; see \cite{Kozak2023}.

We assume that the objective function is quasiconvex and that the constraint functions are strongly quasiconvex. Under the conditions of Lemma \ref{lemma:11}, the latter assumption yields the error bound
\[
 {\rm dist} (x,Q) \leq \frac{2}{\sqrt{\mu}} \sqrt{h^+(x)}, ~~ Q = \{x \in M: \, h(x) \leq 0\}, \; (h^+(x) \equiv \max\{h(x), 0\})
\]
which provides the required control of the distance to the feasible set in terms of constraint violation. Thus, the analysis does not require convexity of the constraint functions.



Let \(f\) be a Lipschitz continuous function with constant \(L_{\Omega}\) on the set \(\Omega = M + {\rm diam} (M) B_{1}(0)\), where \(M\subset \mathbb{R}^{n}\) is a certain convex compact set with
\(\operatorname{diam}(M)\geq 1\). Define the function \(f(x,\alpha)\) as follows:
\[
f(x,\alpha) = \mathbb{E}_{u\sim B_{1}(0)}\left[f(x + \alpha u)\right] = \frac{1}{\operatorname{Vol}(B_{1}(0))}\int_{\| u\| \leq 1}f(x + \alpha u)\mathrm{d}u,
\]
where \(u\) is a random vector uniformly distributed on the unit ball \(B_{1}(0)\) and $\alpha \in [0,1]$.

As it is known, the following results hold  \cite{nesterov2017}:
\begin{enumerate}
 \item $|f(x,\alpha) - f(x)| \leq \alpha L_{\Omega}$ for all $x \in M$.
 
 \item The function \(f(x,\alpha)\) is differentiable.

By differentiating this function and applying Stokes' theorem, we obtain
\begin{align*}
 & \hspace{0.9cm} f^{\prime} (x,\alpha) = \frac{n}{\alpha} \mathbb{E}_{u \sim S^{n - 1}} \left[f(x + \alpha u) u\right] \\
 & \Longrightarrow ~ f^{\prime}(x, \alpha) = -\frac{n}{\alpha} \mathbb{E}_{u \sim S^{n - 1}}\left[f(x - \alpha u) u \right].
\end{align*}
By adding these two equations, we obtain
\begin{align}
 & ~~~ f^{\prime}(x,\alpha) = \frac{n}{2\alpha} \mathbb{E}_{u\sim S^{n - 1}} \left[\left(f(x + \alpha u) - f(x - \alpha u) \right) u\right], \label{eq:11a} \\ 
 & \| f^{\prime} (x,\alpha) \| \leq n L_{\Omega}, ~ \forall ~ x \in M. \label{eq:11b}
\end{align}
\item If the gradient of \(f\) satisfies the Lipschitz condition with constant \(L_{\Omega}^{\prime}\) on the set \(\Omega\), then
\begin{equation}\label{eq:12}
 \| f^{\prime}(x) - f^{\prime} (x, \alpha)\| \leq \alpha L_{\Omega}^{\prime}, ~ \forall ~ x \in M.
\end{equation}
\end{enumerate}

From this representation, it follows that the vector
\[
 g(x,\alpha ,u) = \frac{n}{2 \alpha} \left[f(x + \alpha u) - f(x - \alpha u) \right]u,
\]
is a stochastic gradient of the function \(f(x, \alpha)\) at iteration \(k\). Since, according to the first point, this function is a smooth approximation of \(f\), the vector \(g\) is often referred to in the literature as a stochastic quasi-gradient. 

The following result was proved in \cite{hadjisavvas2026a}.

\begin{theorem} {\rm (\cite[Theorem 4]{hadjisavvas2026a})} \label{theo:03}
 Let \(h\) be a function defined on a convex set \(M\subset \mathbb{R}^{n}\). If \(h\) is strongly quasiconvex with modulus \(\mu >0\), then for every \(x, y \in M\) and every \(z = x + t(y - x)\) with \(0 < t \leq 1\), the following implication holds:
\begin{equation}\label{eq:13}
 h(x) \leq h(z) ~ \Longrightarrow ~ h(z) \leq h(y) - \frac{\mu}{4} (\| y - x\|^{2} - \| z - x\|^{2} ). 
\end{equation}
\end{theorem}

Using the previous theorem, we can prove the following estimate, which will be useful in the sequel.

\begin{lemma} \label{lemma:11}
Assume that the assumptions of Theorem \ref{theo:03} hold.  Let $y\in M$ satisfy $h(y)>0$ and define
\[
 Q := \{x \in M: \, h(x)\leq0\} \neq \emptyset.
\]
Then,
\[
\operatorname{dist}(y,Q)\leq \frac{2\sqrt{h(y)}}{\sqrt{\mu}} .
\]
\end{lemma}

\begin{proof}
 Since $Q$ is nonempty closed and convex, the projection of $y$ onto $Q$ exists and is unique. Denote this projection by $x=P_Q(y)$. Because $h(y)>0$, we have $y\notin Q$. Moreover, by the continuity of $h$, the projection point satisfies $h(x)=0$.

 Take $z=x+t(y-x)$ with $0<t<1$. Since $z$ belongs to the segment connecting a feasible boundary point $x$ and the infeasible point $y$, we have $z\notin Q$, thus $h(z)>0$.

 On the other hand, $h(x)=0<h(z)$. Hence we can apply Theorem \ref{theo:03} with the points $z$ and $x$. Taking $z=x+t(y-x)$, and using the implication in Theorem \ref{theo:03}, we obtain
\begin{align*}
 & h(z)\leq h(y)-\frac{\mu}{4} ( \|y-x\|^2-\|z-x\|^2 ) ~ \Longrightarrow ~ h(z) \leq h(y)-\frac{\mu}{4} (1-t^2)\|y-x\|^2,
\end{align*}
because $\|z-x\|^2 = t^2\|y-x\|^2$. Let $t \to 0$. Then $z\to x$. By the continuity of $h$, $h(z)\to h(x) = 0$. Passing to the limit in the above inequality gives
\begin{align*}
 \frac{\mu}{4} \|y-x\|^2 \leq h(y) ~ \Longrightarrow ~ \|y-x\| \leq \frac{2\sqrt{h(y)}}{ \sqrt{\mu}}.
\end{align*}
Finally, since $x=P_Q(y)$, we have $\|y-x\| = {\rm dist} (y,Q)$.
\end{proof}

\begin{theorem}[Linear error bound under the Slater condition] \label{thm:linear_error_bound}
 Let $M\subset\mathbb{R}^{n}$ be a nonempty closed convex set and let $Q = \{x \in M: \, h(x) \le 0\}$, where $h:M\to\mathbb{R}$ is a continuous strongly quasiconvex function with modulus $\mu>0$. Assume that there exists a Slater point $\bar{x} \in {\rm ri} (M)$ such that
 \[
  h(\bar{x})<0.
 \]
 Then there exists a constant $C>0$ such that
 \[
  {\rm dist} (x, Q) \le C h(x), ~ \forall ~ x \in M \setminus Q.
 \]
 In particular, one may take $C := \frac{2}{\mu \gamma}$ and $\gamma := {\rm dist} (\bar{x}, \partial_M Q) > 0$, where $\partial_M Q$ denotes the boundary of $Q$ relative to $M$.
\end{theorem}

\begin{proof}
 Let $x \in M \setminus Q$. Then $h(x)>0$. Since $h(\bar{x}) < 0$ and $h$ is continuous, $\bar{x} \in {\rm ri}_M (Q)$. Hence, $\gamma := {\rm dist} (\bar{x}, \partial_M Q) > 0$. Consider $z(t)=(1-t) \bar{x}+tx$ with $t \in [0,1]$. By continuity of $h$, there exists a last $t^* \in \, ]0, 1[$ such that $h(z(t^*))=0$. Set $z=z(t^*)$. Then $z\in\partial_M Q$. By strong quasiconvexity,
\[
 0 = h(z) \le h(x)-\frac{\mu}{2} t^*(1-t^*) \|x-\bar{x}\|^2 ~ \Longrightarrow ~ h(x) \ge \frac{\mu}{2} \|x-z\| \|z-\bar{x}\|.
\]
Since $z\in\partial_M Q$, we have $\|z-\bar{x}\|\ge\gamma$. Thus,
\[ 
 h(x) \ge \frac{\mu}{2}\gamma\|x-z\| \ge \frac{\mu}{2} \gamma {\rm dist} (x,Q) ~ \Longrightarrow ~ {\rm dist} (x,Q) \le\frac{2}{\mu\gamma}h(x).
\]
Hence, the assertion follows with $C = \frac{2}{\mu \gamma} > 0$.
\end{proof}

\begin{remark}[On the error bound under quasiconvexity] \label{rem:quasiconvex_slater_error_bound}
 The quasiconvexity of a constraint function together with Slater's condition does not guarantee the square-root error bound used in Lemma \ref{lemma:11} in general. Indeed, let us consider $M=[-1,1]$ and $h(x)=x^3$. Then, $h$ is continuous and strictly increasing, hence quasiconvex, and $Q = \{x\in M: \, h(x) \le 0\} = [-1, 0]$. Slater's condition is satisfied, for instance, at $\bar x=-1/2$, since $h(\bar x)=-1/8<0$. For $x>0$, however, $\operatorname{dist}(x,Q)=x$ and $h^+(x)=x^3$. Consequently, an estimate of the form
\[
\operatorname{dist}(x,Q)\le C\sqrt{h^+(x)}
\]
would imply $x\le Cx^{3/2}$, thus $x^{-1/2}\le C$, which is impossible as $x\downarrow0$. Thus, quasiconvexity and Slater's condition alone are insufficient to ensure the error bound required in
our analysis.

Under the strong quasiconvexity assumption of Lemma \ref{lemma:11}, the required estimate does hold. More precisely, if $Q = \{x\in M: \, h(x)\le0\}$ and $h$ has a strictly feasible point, then
\[
 {\rm dist} (y,Q) \le \frac{2}{\sqrt{\mu}} \sqrt{h^+(y)}.
\]
This estimate provides the quantitative relation between constraint
violation and distance to the feasible set needed in the convergence
analysis.
\end{remark}

Now, let us consider the constrained optimization problem
\begin{equation}\label{eq:14}
 \min f_{0}(x),\quad x \in M, \quad f_{i}(x) \leq 0,\; i = 1,2,\ldots ,m. 
\end{equation}

For which we suppose that 
$Q=\{ x \in M: f_i(x) \leq 0, i=1,2,..., m\} \neq \emptyset$. Define the function
\[
\tilde{f} (x) := \max_{1 \leq i \leq m} f_{i}(x).
\]


In this section, we present a Feasibility--Optimization Switching
Stochastic Finite-Difference Method for Problem~\eqref{eq:14} and
establish the corresponding convergence guarantees.

\subsection{Algorithmic Scheme}

The proposed method utilizes an adaptive switching mechanism that alternates between minimizing the objective function and reducing constraint violations based on an asymptotically vanishing feasibility threshold.

\begin{algorithm}[H]
\caption{Feasibility--Opt Switching Stochastic Finite-Difference Method}
\label{alg:stochastic_switching}
\begin{algorithmic}[1]

\Require
\Statex $\bullet$ $x^0 \in M$, $\rho_k = \frac{1}{\sqrt{k+0.5}}$, $\alpha_k = \frac{1}{k+1}$, $R = \operatorname{diam}(M) = \max_{x,y \in M} \|x-y\|$,
 common Lipschitz constant on $\Omega$, $L_\Omega = \max_{0 \leq i \leq m} L_i > 0$, where $L_i$ is the Lipschitz constant of $f_i$ on $\Omega$, and a positive lower bound on gradient norms:  $\gamma = \min_{0\leq i\leq m} \min_{x\in M} \|f_i'(x)\| > 0$, 
\Statex $\bullet$ Unified scale parameter:
\[
\widetilde R \geq \max\left\{ R,\; n L_\Omega,\; \frac{L_\Omega}{2\gamma}\left(n^2 L_\Omega^2 + \frac{2R^2}{\ln 3} + 1\right) \right\}
\]
\For{$k=0,1,2,\ldots$}

    \State Compute the maximum constraint violation:
    \[
    \widetilde f(x^k) = \max_{1\leq i\leq m} f_i(x^k).
    \]

    \If{$\widetilde f(x^k) \leq \rho_k \widetilde R$}

        \State Sample a direction uniformly from the unit sphere: $u_k \sim \operatorname{Unif}(S^{n-1})$.
        \State Set search direction for objective minimization:
        \[
        p^k = g_0(x^k, \alpha_k, u_k).
        \]

    \Else

        \State Identify active maximal-violation constraints:
        \[
        I(x^k) = \left\{ i \in \{1,\ldots,m\} : f_i(x^k) = \widetilde f(x^k) \right\}.
        \]
        \State Select an index $\xi_k$ uniformly at random from $I(x^k)$.
        \State Sample a direction uniformly from the unit sphere: $u_k \sim \operatorname{Unif}(S^{n-1})$.
        \State Set search direction for constraint reduction:
        \[
        p^k = g_{\xi_k}(x^k, \alpha_k, u_k).
        \]

    \EndIf

    \State Update the iterate via projection onto $M$:
    \[
    x^{k+1} = \Pi_M (x^k - \rho_k p^k).
    \]

\EndFor

\Ensure A sequence of iterates $\{x^k\}_{k} \subset M$.
\end{algorithmic}
\end{algorithm}

\subsection{Convergence Result}

The main theoretical outcome establishing the rate of convergence for Algorithm~\ref{alg:stochastic_switching} is formalized below.

\begin{theorem}\label{theo:convergence}
 Consider the constrained optimization problem on $M \subset \mathbb{R}^n$. Assume the following conditions hold:
\begin{enumerate}
    \item $M \subset \mathbb{R}^n$ is a nonempty, convex, and compact set with diameter $R= \operatorname{diam}(M)$. Define the enlarged domain 
    \[
    \Omega = M + R B_1(0).
    \]
    
    \item The objective function $f_0$ is quasiconvex on $\Omega$, while each constraint function $f_i$ ($i = 1,\ldots,m$) is $\mu$-strongly quasiconvex on $\Omega$ for some $\mu > 0$.
    
    \item Each function $f_i$, $i=0,\ldots,m$, has a Lipschitz-continuous gradient on $\Omega$. Consequently, each $f_i$ is Lipschitz continuous on $\Omega$. We denote by $L_\Omega>0$ a common Lipschitz constant for the functions $f_0, \ldots, f_m$ on $\Omega$.
    
    \item Gradients do not vanish on $M$, admitting a strictly positive lower bound:
    \[
     \gamma := \min_{0 \leq i \leq m} \min_{x \in M} \|f'_i(x)\| > 0. 
    \]
\end{enumerate}

Let the algorithm parameters be chosen as
\[
 \rho_k = \frac{1}{\sqrt{k+0.5}}, \qquad \alpha_k = \frac{1}{k+1}, 
\]

Then, for the sequence $\{x^k\}_{k}$ generated by Algorithm~\ref{alg:stochastic_switching}, there exists a constant $C>0$ such that, for every integer $k\geq 3$, one can find an index $i'\in\{0,1,\ldots,k\}$ satisfying
\[
 \lvert f_0(x^{i'}) - f_0(x^*) \rvert \leq \frac{A}{k^{1/4}} \qquad  {\rm and} \qquad d(x^{i'}, Q) \leq \frac{A}{k^{1/4}},
\]
where $x^*$ is an optimal solution to the constrained problem, $Q = \{x \in M: \, f_i(x) \leq 0, \, i=1,\ldots,m\}$ is the feasible set, and $d(x, Q) = \inf_{y \in Q} \|x - y\|$ represents the Euclidean distance to $Q$.
\end{theorem}

\begin{proof}
 Let \(L'_\Omega\) denote a common Lipschitz constant for the gradients of \(f_0, f_1, \ldots, f_m\) on \(\Omega\). For every \(k\), define
\[
 I(x^k) = \{i \in \{1,\ldots, m\}: \, f_i(x^k) = \widetilde f(x^k)\}.
\]

First consider the case
\[
 \widetilde f(x^k) > \rho_k \widetilde R.
\]

In this case, the algorithm selects $\xi_k$ uniformly from $I(x^k)$, while $u_k$ is uniformly distributed on $S^{n-1}$. By the representation of the stochastic finite-difference estimator,
its conditional mathematical expectation with respect to $x^k$ is the corresponding smoothed gradient. Hence
\[
 \mathbb{E}[p^k \mid x^k] = \frac{1}{ \lvert I(x^k) \rvert} \sum_{i \in I(x^k)} f_i'(x^k, \alpha_k).
\]

For every \(i\in I(x^k)\), $f_i(x^k)=\widetilde f(x^k)$. Then, by the smoothing estimate, 
\[
\begin{aligned}
 \langle f_i'(x^k,\alpha_k),x^k-x^* \rangle & \geq \langle f_i'(x^k),x^k-x^* \rangle - L'_\Omega \operatorname{diam} (M) \alpha_k.
\end{aligned}
\]

Since \(x^*\in Q\), we have $ f_i(x^*)\leq0$. Moreover, \(f_i(x^k)=\widetilde f(x^k)\). Hence, by the geometric inequality,
\[
\begin{aligned}
 \langle f_i'(x^k), x^k-x^* \rangle & \geq \frac{ \gamma}{L_\Omega} ( f_i(x^k)-f_i (x^*)) \geq \frac{\gamma}{L_\Omega} \widetilde f(x^k).
\end{aligned}
\]

Consequently,
\[
 \langle f_i'(x^k,\alpha_k),x^k-x^* \rangle \geq \frac{\gamma}{L_\Omega} \widetilde f(x^k) - L'_\Omega \operatorname{diam}(M)\alpha_k.
\]

This estimate holds for every \(i\in I(x^k)\). Then, after averaging over all \(i\in I(x^k)\), we obtain
\begin{align*}
 \langle \mathbb{E}[p^k \mid x^k], x^k-x^* \rangle & \geq \frac{\gamma}{L_\Omega} \widetilde f(x^k) - L'_\Omega {\rm diam} (M) \alpha_k  > \frac{\gamma \widetilde R}{L_\Omega} \rho_k - L'_\Omega {\rm diam} (M) \alpha_k,
\end{align*}
where the strict inequality follows from $\widetilde f(x^k)>\rho_k\widetilde R$.

Consider the projected iteration
\[
 x^{i+1} = \Pi_M(x^i-\rho_ip^i).
\]
Since \(x^*\in M\) and $M$ is closed and convex, the nonexpansiveness of the projection gives
\[
\begin{aligned}
 \|x^{i+1} - x^*\|^2 & \leq \|x^i - x^* - \rho_i p^i\|^2 = \|x^i - x^*\|^2 - 2\rho_i \langle p^i,x^i-x^*\rangle + \rho_i^2 \|p^i\|^2.
\end{aligned}
\]

Since $\|p^i\|\leq nL_\Omega$,
we have
\[
 \|x^{i+1}-x^*\|^2 \leq \|x^i-x^*\|^2 - 2\rho_i \langle p^i, x^i-x^*\rangle + n^2L_\Omega^2 \rho_i^2.
\]

Taking the conditional mathematical expectation with respect to \(x^i\), we obtain
\[
\begin{aligned}
 \mathbb{E} \left[ \|x^{i+1}-x^*\|^2 \mid x^i \right]
 & \leq \|x^i-x^*\|^2 - 2\rho_i \left\langle \mathbb{E} [p^i\mid x^i], x^i-x^* \right\rangle + n^2 L_\Omega^2\rho_i^2.
\end{aligned}
\]

Put $r_i=\|x^i-x^*\|$. Let $k'=\lfloor k/3\rfloor$ and define
\[
 I_k = \{i \in\{k',\ldots,k\}:\widetilde f(x^i) \leq \rho_i \widetilde R\}.
\]

Thus, \(I_k\) is the set of iterations at which the algorithm operates in the objective-minimization mode.

First, we prove that \(I_k \neq \emptyset\). Suppose, to the contrary, that \(I_k = \emptyset\). Then for all \(i \in \{k', \dots, k\}\), we have \(\widetilde f(x^i) > \rho_i \widetilde R\).

Summing the recurrence inequality for \(\mathbb{E}[r_{i+1}^2 \mid x^i]\) from \(i=k'\) to \(i=k\) and taking total expectations gives:
\[
\mathbb{E}[r_{k+1}^2] \leq \mathbb{E}[r_{k'}^2] - \left(\frac{2\gamma \widetilde R}{L_\Omega} - n^2 L_\Omega^2\right) \sum_{i=k'}^k \rho_i^2 + 2L'_\Omega R \sum_{i=k'}^k \rho_i \alpha_i.
\]
By the choice of \(\widetilde R\), we have
\[
\frac{2\gamma \widetilde R}{L_\Omega} - n^2 L_\Omega^2 \geq \frac{2R^2}{\ln 3} + 1.
\]
Furthermore, \(\mathbb{E}[r_{k'}^2] \le R^2\), and for all sufficiently large \(k\),
\[
\sum_{i=k'}^k \rho_i^2 \geq \frac{\ln 3}{2},
\qquad
\sum_{i=k'}^k \rho_i \alpha_i \leq \frac{C_0}{\sqrt{k}}.
\]
Substituting these estimates yields
\[
\begin{aligned}
0 \leq \mathbb{E}[r_{k+1}^2]
& \leq R^2 - \left(\frac{2R^2}{\ln 3} + 1\right)\frac{\ln 3}{2} + \frac{C}{\sqrt{k}} \\
& = -\frac{\ln 3}{2} + \frac{C}{\sqrt{k}}.
\end{aligned}
\]
Since \(\ln 3/2 > 0\), the right-hand side is strictly negative for all sufficiently large \(k\), which is a contradiction.

Thus, \(I_k \neq \emptyset\). Now, for \(i\in I_k\), define
\[
v_i
=
\left\langle
f_0'(x^i,\alpha_i),x^i-x^*
\right\rangle.
\]

For \(i\notin I_k\), define
\[
\widetilde v_i
=
\left\langle
\mathbb{E}[p^i\mid x^i],
x^i-x^*
\right\rangle.
\]

Summing the conditional expectations from \(i=k'\) to \(i=k\) and taking total expectations gives
\[
\begin{aligned}
\mathbb{E} [r_{k+1}^2] & + 2\sum_{i\in I_k}\rho_i \mathbb{E} [v_i] + 2\sum_{i\notin I_k}\rho_i \mathbb{E} [\widetilde v_i] \leq \mathbb{E}[r_{k'}^2] + n^2L_\Omega^2 \sum_{i=k'}^k \rho_i^2.
\end{aligned}
\]

Suppose, to the contrary, that $v_i\geq\widetilde R^2\rho_i$ for all $i\in I_k$. Using the definition of \(\widetilde R\), the lower bound for \(\widetilde v_i\), and the estimates for the constraint-reduction iterations, we obtain:
\[
1
\geq
\sum_{i=k'}^k\frac{1}{i+0.5}
-
2L'_\Omega\operatorname{diam}(M)
\sum_{i=k'}^k\rho_i\alpha_i.
\]

Since $\sum_{i=k'}^k\frac{1}{i+0.5} \geq \ln3$, and $\sum_{i=k'}^k \rho_i\alpha_i \leq \frac{C_0}{\sqrt{k}}$, we conclude that $1 \geq \ln3-\frac{C}{\sqrt{k}}$, where \(C>0\) is independent of \(k\). Since \(\ln3>1\), the right-hand side is greater than \(1\) for all sufficiently large \(k\), which is a contradiction. Therefore, there exists an index \(i'\in I_k\) such that $\langle f_0'(x^{i'}, \alpha_{i'}), x^{i'}-x^* \rangle < \widetilde R^2 \rho_{i'}$.

First suppose that $x^{i'}\in Q$. Using the smoothing estimate and the geometric inequality for the strongly quasiconvex objective function \(f_0\), we obtain
\[
\begin{aligned}
 \widetilde R^2 \rho_{i'} & > \langle f_0' (x^{i'}, \alpha_{i'}), x^{i'}-x^* \rangle \geq \frac{\gamma}{L_\Omega} (f_0(x^{i'}) -f_0(x^*)) - L'_\Omega {\rm diam} (M) \alpha_{i'}.
\end{aligned}
\]

Hence,
\[
 f_0(x^{i'}) - f_0(x^*) \leq C_1 \rho_{i'} + C_2 \alpha_{i'}.
\]

Since $i'\geq k'$ and $k' > \frac{k}{3} - 1$, there exists $C_3 > 0$ such that $\rho_{i'} \leq \frac{C_3}{\sqrt{k}}$. Then,
\[
 f_0 (x^{i'}) - f_0(x^*) \leq \frac{C_4}{\sqrt{k}}.
\]

Now consider the case $x^{i'} \notin Q$. Let $y^{i'} = \Pi_Q(x^{i'})$. By Lemma \ref{lemma:11},
\[
 \|x^{i'} - y^{i'}\| \leq \frac{2}{\sqrt{\mu}} \sqrt{ \widetilde f^+(x^{i'})}.
\]

By expanding the inner product difference, we have
\[
\begin{aligned}
 & ~~\,~ \langle f_0'(x^{i'}, \alpha_{i'}), x^{i'}-x^* \rangle - \langle f_0'(y^{i'}, \alpha_{i'}), y^{i'} -x^* \rangle
 \\
 & = \langle f_0'(x^{i'}, \alpha_{i'}) - f_0' (y^{i'}, \alpha_{i'}), x^{i'}-x^* \rangle + \langle f_0'(y^{i'}, \alpha_{i'}), x^{i'}-y^{i'} \rangle.
\end{aligned}
\]
Using the Lipschitz continuity of the gradient of \(f_0\) (with constant \(L'_\Omega\)) for the first term and the boundedness of the smoothed gradient \(\|f_0'(y^{i'},\alpha_{i'})\| \leq L_\Omega\) for the second term, we obtain
\[
\begin{aligned}
 \langle f_0'(x^{i'}, \alpha_{i'}), x^{i'}-x^* \rangle & \geq \langle f_0' (y^{i'}, \alpha_{i'}), y^{i'} - x^* \rangle  - (L'_\Omega {\rm diam} (M) + L_\Omega) \|x^{i'}-y^{i'}\| \\
 & \geq \langle f_0' (y^{i'}, \alpha_{i'}), y^{i'} - x^* \rangle -C_5 \sqrt{\widetilde f^+(x^{i'})}.
\end{aligned}
\]

Since \(y^{i'}\in Q\), we may apply the estimate obtained above for the objective function. Hence,
\[
\begin{aligned}
 \langle f_0'(x^{i'}, \alpha_{i'}), x^{i'}-x^* \rangle & \geq \frac{\gamma}{L_\Omega} ( f_0(x^{i'}) - f_0(x^*)) - C_6 \sqrt{\widetilde f^+(x^{i'})}.
\end{aligned}
\]

On the other hand, by the choice of \(i'\), we have
$\langle f_0' (x^{i'}, \alpha_{i'}), x^{i'}-x^* \rangle < \widetilde R^2 \rho_{i'}$. Then,
\[
 f_0(x^{i'}) - f_0(x^*) \leq C_7 \sqrt{\widetilde f^+(x^{i'})} + C_8 \rho_{i'}.
\]

Since \(i'\in I_k\), the algorithm is in the objective-minimization mode at iteration $i'$. Hence, by the switching rule, $\widetilde f(x^{i'}) \leq \rho_{i'}\widetilde R$. Consequently,
\[
 \widetilde f^+(x^{i'}) \leq \widetilde R\rho_{i'}.
\]
Since \(i' \geq k'\), there exists \(C_9>0\) such that $\rho_{i'} \leq \frac{C_9}{\sqrt{k}}$. Then, $\sqrt{ \widetilde f^+(x^{i'})} \leq \frac{C_{10}}{k^{1/4}}$. It follows that there exists $A_1>0$, independent of $k$, such that
\[
 f_0(x^{i'})-f_0(x^*) \leq \frac{A_1}{k^{1/4}}.
\]
Moreover, for every \(i=1,\ldots,m\), we have $f_i(x^{i'}) \leq \widetilde f(x^{i'}) \leq \frac{C_{11}}{\sqrt{k}}$. Applying Lemma~\ref{lemma:11}, we obtain
\[
\begin{aligned}
 {\rm dist} (x^{i'},Q) &\leq \frac{2}{\sqrt{ \mu}} \sqrt{\widetilde f^+(x^{i'})} \leq \frac{A_2}{k^{1/4}}.
\end{aligned}
\]

It remains to obtain a lower bound for the objective error.

Since $y^{i'}\in Q$ and \(x^*\) is an optimal solution, we have $f_0(y^{i'})\geq f_0(x^*)$. Since $f_0$ is Lipschitz continuous on \(\Omega\), we have
\begin{align*}
 & \hspace{0.85cm} f_0 (x^{i'}) - f_0 (y^{i'}) \geq - L_\Omega \| x^{i'} - y^{i'}\| \\
 & \Longrightarrow ~ f_0(x^{i'})-f_0(x^*)  \geq - L_\Omega {\rm dist} (x^{i'},Q) \geq -\frac{A_3}{k^{1/4}}.
\end{align*}

Combining this estimate with the previously obtained upper bound, we conclude that
\[
 \left| f_0(x^{i'})-f_0(x^*) \right| \leq \frac{A}{k^{1/4}},
\]
for some constant \(A>0\) independent of $k$.

At the same time, $\operatorname{dist}(x^{i'},Q) \leq \frac{A}{k^{1/4}}$. Thus, for every \(k\geq3\), there exists \(i'\in\{0,\ldots,k\}\) for which both estimates hold. The proof is complete.
\end{proof}

\begin{remark}
If the Slater condition is additionally imposed, then the standard
linear error bound
\[
\operatorname{dist}(x,Q)\leq C_S h^{+}(x)
\]
allows the constraint-violation estimate in
Theorem~\ref{theo:convergence} to be converted into a
distance-to-feasibility estimate of the same order.

Moreover, under the additional estimates used in the proof of
Theorem~\ref{theo:convergence}, the objective-function error satisfies
\[
|f_0(x^{i'})-f_0(x^*)|
\leq \frac{C_S'}{\sqrt{k}},
\]
where $C_S'>0$ is independent of $k$. Hence, to achieve
\[
|f_0(x^{i'})-f_0(x^*)|\leq\varepsilon,
\]
it is sufficient to take
\[
k=O(\varepsilon^{-2}).
\]

For comparison, Nguyen and Balasubramanian \cite{NguyenBalasubramanian2022} establish an oracle complexity of order $O(\varepsilon^{-2})$ in the convex setting and $O(\varepsilon^{-3})$ in the nonconvex setting, with respect to their approximate-KKT criterion. Thus, with respect to the accuracy parameter, the present Slater-based objective-value estimate has the same order $O(\varepsilon^{-2})$ as their convex result.

This comparison is only at the level of the dependence on the accuracy parameter. The two results are obtained under different assumptions and convergence criteria. In particular, our objective function is assumed to be quasiconvex, whereas in \cite{NguyenBalasubramanian2022} consider a different class of functionally constrained stochastic optimization problems. Moreover, their approximate-KKT guarantee does not directly imply the objective-function value estimate established in our theorem.
\end{remark}

\begin{remark} \label{rem:interpretation_convergence}
 \begin{itemize} 
  \item[$(i)$] Under the assumptions of Theorem \ref{theo:convergence}, the asymptotic optimality and convergence to the feasible set $Q$ can be equivalently characterized in terms of the objective distance and set distance as:
\[
 \liminf_{k \to + \infty} \mathbb{E} \left[ \left| f_0(x^k) - f_0(x^*) \right|\right]=0, ~~ \text{and} ~~ \liminf_{k \to + \infty} \mathbb{E} \left[{\rm dist} (x^k, Q) \right] = 0.
\]
 
 \item[$(ii)$] A similar algorithm was proposed by Nesterov in \cite{nesterov2010} for the convex case, that is, when all involved functions are convex and their subgradients are known (\cite[Theorems 3.2.3, p. 177]{nesterov2010}), in which case a convergence rate of $\frac{1}{\sqrt{k}}$ was obtained. 
 \end{itemize} 
\end{remark}

\section{The Penalty Function Method}\label{sec:05}

The presentation of this method is based on problem \eqref{eq:14}. We assume that all functions \(f_i\), with \(i=1,\ldots,m\), are strongly quasiconvex and locally Lipschitz continuous on \(\mathbb{R}^n\), that the function \(f_0\) is quasiconvex and  regular, and that the feasible set \(M \subset \mathbb{R}^n\) is convex and compact.

We also assume that the minimal value
\[
 f_{0}^{*} := \min_{x\in Q} f_{0}(x),\quad Q := \{x \in M: ~ f_{i}(x) \leq 0, ~ i = 1,\ldots ,m\},
\]
is known. The goal is to generate a sequence \(\{x^{k}\}_{k} \subset M\) that converges to the set
\[
 Q^{*} := \{x\in Q:f_{0}(x) = f_{0}^{*}\}
\]
of minimizers. Consider the auxiliary function,
\[
 \tilde{f}_{0}(x) := \max \{f_{0}(x) - f_{0}^{*}, f_{1}(x), \ldots, f_{m}(x)\}.
\]
Clearly, \(\tilde{f}_{0}\) is quasiconvex and \(Q^{*}\) is the set of minimizers of \(\tilde{f}_{0}\) on \(Q\).

The algorithm performs a sliding minimization over the family of functions
\[
 F(x,C) := \tilde{f}_{0}(x) + C\sum_{i = 1}^{m} \left[\max \{0,f_{i}(x)\}\right]^{2}, \quad C > 0,
\]
on \(M\), incorporating a quadratic penalty. 

Let the sequence of iterates be generated by
\[
x^{k + 1} = \Pi_{M} \big(x^{k} - \alpha_{k}a^{k} \big), ~ \forall ~ k = 0,1,\ldots
\]
where  $a^{k}\in \partial F(x^{k},C_{k})$ and, more explicitly,
\[
 a^{k} = a_{0}^{k} + 2C_{k}\sum_{i = 1}^{m}\max \{0,f_{i}(x^{k})\} a_{i}^{k}, ~~ a_i^0 \in \partial \tilde{f}_0(x^k), ~~ a_{i}^{k}\in \partial f_{i}(x^{k}),\; i = 1,\ldots ,m.
\]

We impose the following conditions on the parameters
\[
 \alpha_{k} > 0, ~~ \sum_{k = 0}^{\infty} \alpha_{k} = + \infty, ~~ C_{k} \to + \infty, ~~ \sum_{k = 0}^{\infty} \alpha_{k}^{2} C_{k}^{2} < + \infty, ~~ \sum_{k = 0}^{\infty} \frac{ \alpha_{k}}{C_{k}^{1/3}} < + \infty.
\]

For example, the previous assumptions hold when $\alpha_{k} = \frac{1}{k}$ and $C_{k} = k^{1 / 4}$.

We have the following convergence result.

\begin{theorem}\label{theo:05}
 Let us consider the assumptions mentioned above on the functions $f_{i}$, with $i = 0,1 \ldots ,m$, and the set $M$. Assume that 
 \[
  0 \notin \partial f_{i}(x), ~ \forall ~ x \in M, ~ \forall ~ i = 0,1,\ldots, m.
 \]
 Then the sequence $\{x^{k}\}_{k}$ generated by the above sliding minimization algorithm converges to the set $Q^{*}$.
\end{theorem}

\begin{proof}
 Let \(x^{*} \in Q^{*}\). Then the projection update yields
\begin{equation}\label{eq:16}
 \| x^{k + 1} - x^{*}\|^{2} \leq \| x^{k} - x^{*} \|^{2} - 2\alpha_{k}\langle a^{k}, x^{k} - x^{*} \rangle + \alpha_{k}^{2}\| a^{k}\|^{2}. 
\end{equation}

We first estimate the middle term using Lemma \ref{lemma:04}. Let \(i \in \{1, \ldots, m\}\) be with \(f_{i}(x^{k}) > 0\). Then by inequality \eqref{eq:03},
\[
 f_{i}(x^{k}) - f_{i}(x^{*}) \leq L_{\Omega} \left\langle \frac{a_{i}^{k}}{\|a_{i}^{k}\|},x^{k} - x^{*} \right\rangle \leq \frac{L_{\Omega}}{\tilde{\gamma}} \langle a_{i}^{k},x^{k} - x^{*} \rangle ,
\]
where $\tilde{\gamma} := \min_{i} \, \min_{x\in M} \,  \operatorname{dist}(0,\partial f_{i}(x)) > 0$. Thus,
\[
 2 C_{k} \sum_{i = 1}^{m} \max\{0,f_{i} (x^{k})\} \langle a_{i}^{k},x^{k} - x^{*} \rangle \geq \frac{2 \tilde{\gamma}C_{k}}{L_{\Omega}} \sum_{i = 1}^{m} \max \{0, f_{i} (x^{k})\} (f_{i}(x^{k}) - f_{i}(x^{*})).
\]

Since \(\tilde{f}_{0}(x^{*}) = 0\) and \(\tilde{f}_{0}(x^{k}) > 0\), we have
\[
 \langle a_0^k,x^k - x^*\rangle \geq \frac{ \tilde{\gamma}}{L_{\Omega}} (\tilde{f}_0 (x^k) - \tilde{f}_0 (x^*)).
\]

Combining these, we obtain
\[
 \langle a^{k},x^{k} - x^{*}\rangle \geq \frac{2 \tilde{\gamma}}{L_\Omega} (F(x^{k},C_{k}) - F(x^{*}, C_{k})).
\]

Next, we estimate the difference \(F(x^{k},C_{k}) - F(x^{*}, C_{k})\) using Lemma \ref{lemma:11}. Note that \(F(x,C_{k}) = \tilde{f}_{0}(x)\) for all \(x\in Q\). Hence,
\[
 F^{*} (C_{k}) := \inf_{x \in M} F(x,C_{k}) \leq \inf_{x\in Q} F(x,C_{k}) = \tilde{f}_{0}^{*}.
\]

For \(x\in M\setminus Q\), we have for any \(y\in Q\) that
\[
 \tilde{f}_{0} (x) \geq \tilde{f}_{0}(y) - L\| x - y\| \geq \tilde{f}_{0}^{*} -  L ({\rm dist} (x, Q)).
\]

Moreover, by Lemma \ref{lemma:11},
\[
 C_{k} \sum_{i = 1}^{m}[\max \{0,f_{i}(x)\} ]^{2} \geq \frac{C_{k}\mu^{2}}{8} (\operatorname{dist}(x,Q))^{4}.
\]

Combining the above estimates, there exists a constant \(A > 0\) such that
\[
\tilde{f}_{0}^{*} - F(x,C_{k}) \leq \frac{A}{C_{k}^{1 / 3}}.
\]

Since the subgradients $a^{k}$ are bounded on the compact set $M$, there exists $B > 0$ such that $\| a^{k}\| \leq B C_{k}$. Substituting these bounds into \eqref{eq:16} gives
\begin{equation}\label{eq:17}
 \| x^{k + 1} - x^{*}\|^{2} \leq \| x^{k} - x^{*} \|^{2} + \alpha_{k} \frac{A \tilde{\gamma}}{L_{\Omega}C_{k}^{1 / 3}} + B^{2} C_{k}^{2} \alpha_{k}^{2}. 
\end{equation}

Since, by assumptions, the series $\sum_{k = 0}^{\infty} (\frac{\alpha_{k}}{C_{k}^{1 / 3}} + \alpha_{k}^{2} C_{k}^{2}) < + \infty$, thus ${\rm dist} (x^{k}, Q^{*})$ is convergent. Furthermore, from \eqref{eq:17} there exists a subsequence $\{x^{k_{j}}\}_{j}$ such that
\[
 F(x^{k_{j}}, C_{k_{j}}) = \tilde{f}_{0}(x^{k_{j}}) + C_{k_{j}} \sum_{i = 1}^{m} (\max \{0,f_{i}(x^{k_{j}})\})^{2} \, \to \, \tilde{f}_{0}^{*}.
\]

Since this subsequence lies in the compact set $M$, we may assume it converges to some $x^{**} \in M$. Then, $\sum_{i = 1}^{m} [\max \{0, f_{i} (x^{**})\} ]^{2} = 0$ and $\tilde{f}_{0} (x^{**}) = \tilde{ f}_{0}^{*}$, so that \(x^{**} \in Q^{*}\). Hence, ${\rm dist} (x^{k_{j}}, Q^{*}) \to 0$. Since the sequence ${\rm dist} (x^{k}, Q^{*})$ is convergent, it follows that ${\rm dist} (x^{k}, Q^{*}) \to 0$, which completes the proof.
\end{proof}

\begin{remark}
 It follows from the proof of Theorem \ref{theo:05} that, if the objective function $f_0$ is convex, then the conclusion remains valid even when its minimal value is unknown. In this case, it suffices to set $\tilde{f}_0(x)=f_0(x).$
Then, the inequality
\[
 \langle a_k, x^k - x^* \rangle \ge \frac{2 \tilde{\gamma}}{L_{\Omega}} (\tilde{f}_0 (x^k) - \tilde{f}_0 (x^*) )
\]
can be replaced by the fundamental inequality of convex functions
\[
 \langle a_k, x^k - x^* \rangle \ge {f}_0(x^k) - {f}_0(x^*),
\]  
which holds automatically. Then we just repeat the proof of the Theorem \ref{theo:05}.
\end{remark}

\begin{remark}
The penalty function method replaces a constrained problem by an unconstrained one, adding to the objective a term that penalizes constraint violations. It is simple to implement, avoids projections onto the feasible set, and can be combined with standard (stochastic) gradient methods. As the penalty parameter grows, solutions of the penalized problem approach those of the original one, although very large parameters may cause ill-conditioning and slow convergence. In contrast to \cite[Theorem 5.5, p.~256]{Sukharev1986}, where convergence to the set of minimizers is proved under convexity of all functions, Theorem~\ref{theo:05} only requires quasiconvexity, together with the condition $0 \notin \partial f_i(x)$, $i=0,1,\ldots,m$, which is not assumed in \cite{Sukharev1986}.
\end{remark}

We conclude this paper with a numerical illustration of the finite difference stochastic method (Algorithm~1).

\begin{example}\label{ex:final}
Consider the constrained optimization problem
\begin{equation}\label{eq:problem}
    f_0(x_1,x_2) = x_1+x_2+0.01(x_1^2+x_2^2) \to \min
\end{equation}
subject to
\begin{equation}\label{eq:constraint}
    f_1(x_1,x_2) = 0.8x_1^2+0.4x_2^2 - \frac{1}{0.8x_1^2+0.4x_2^2+10} - 1 \leq 0,
    \quad (x_1,x_2)\in M = B_{0.45}(1,1).
\end{equation}
The objective $f_0$ is strongly convex (hence quasiconvex), and $f_1$ is strongly quasiconvex by Remark~1(ii).

With diameter $R = \operatorname{diam}(M) = 0.9$, let $\Omega = M + RB_1(0) = B_{1.35}(1,1)$ ($n=2$). 
The gradient bounds over $\Omega$ are $\gamma_0 \approx 1.4335$, $L_0 \approx 1.4695$, and for the constraint $f_1$, numerical estimates yield $\gamma_1 \approx 1.1552$, $L_1 \approx 3.8754$. 
Setting $\gamma = \min\{\gamma_0,\gamma_1\} \approx 1.1552$ and $L_\Omega = \max\{L_0,L_1\} \approx 3.8753$, we choose the unified scale parameter $\widetilde{R} = 105$.

Parameters are $\rho_k = 1/\sqrt{k+0.5}$, $\alpha_k = 1/(k+1)$, with initial point $x^0 = (1.3, 0.7)$ ($f_0(x^0) = 2.0218$, $f_1(x^0) \approx 0.4614 > 0$, thus infeasible).

The global minimizer is $x^* \approx (0.6818, 0.6818)$ with $f_0(x^*) \approx 1.3729$ and $f_1(x^*) \approx -0.5369 < 0$. The feasible set is $Q = M \cap \{x : 2x_1^2+x_2^2 \leq 2.7254\}$.

Numerical Results: After $5000$ iterations of Algorithm~1 implemented in Python (using a fixed random seed), the results are:
\begin{equation*}
    x^{5000} \approx (0.6813, 0.6823), \quad f_0(x^{5000}) \approx 1.372902, \quad f_1(x^{5000}) \approx -0.5372 < 0.
\end{equation*}
The errors are $\|x^{5000}-x^*\| \approx 0.00071$ and $|f_0(x^{5000})-f_0(x^*)| \approx 8.04 \times 10^{-7}$. The objective-minimization phase remained active throughout, while constraint-reduction was inactive.
\end{example}

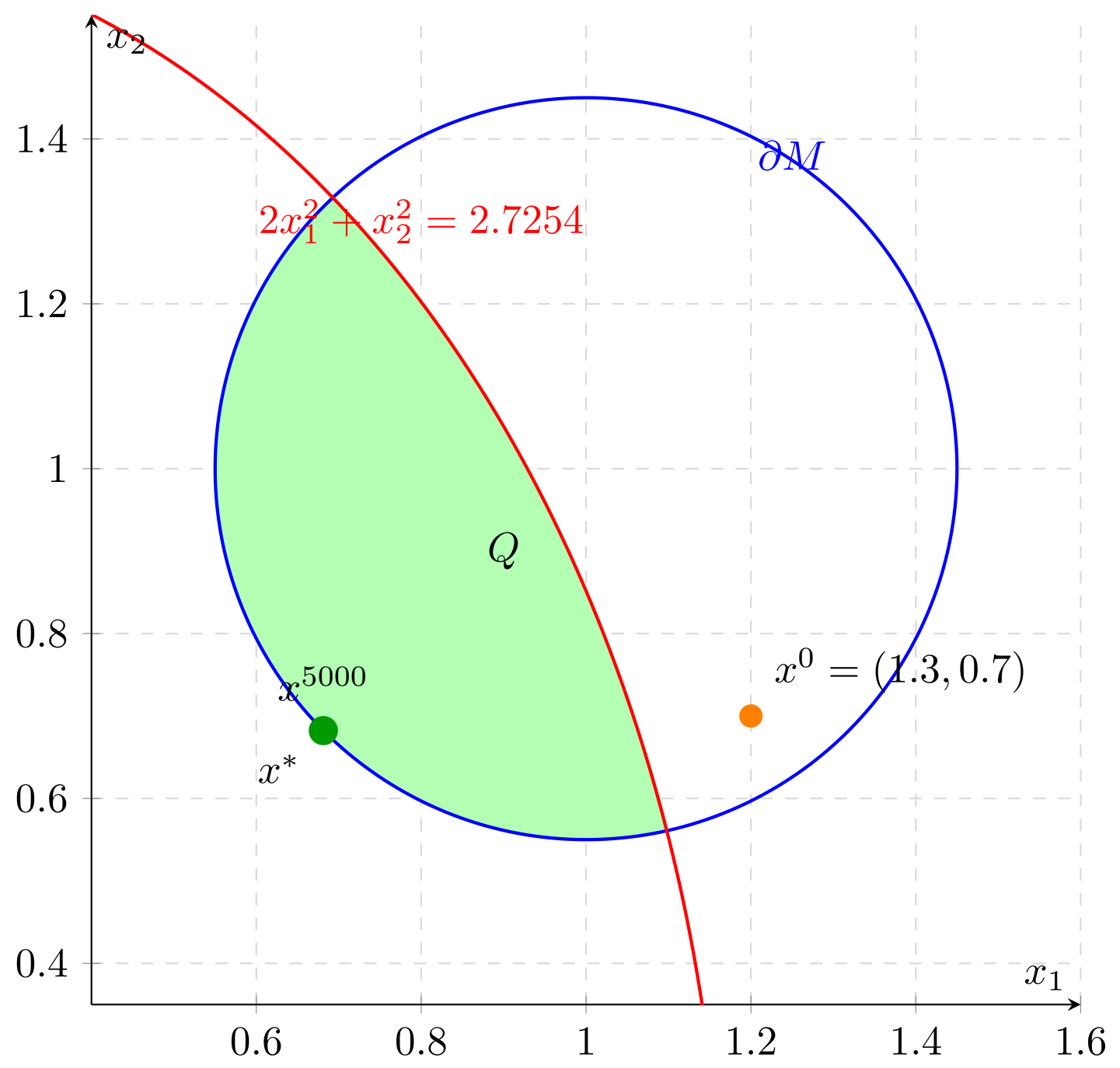
\begin{figure}[H]
\centering
\begin{tikzpicture}
\begin{axis}[
    axis lines=middle,
    xlabel={$x_1$},
    ylabel={$x_2$},
    xmin=0.4, xmax=1.6,
    ymin=0.4, ymax=1.5,
    grid=major,
    grid style={dashed,gray!30},
    width=10cm,
    height=10cm,
    axis equal=true,
    clip=true
]

\begin{scope}
  \clip
  plot[domain=0.4:1.167, samples=200, variable=\x]
    (axis cs:\x,{sqrt(max(0,2.7254-2*\x*\x))})
  -- plot[domain=1.167:0.4, samples=200, variable=\x]
    (axis cs:\x,{-sqrt(max(0,2.7254-2*\x*\x))})
  -- cycle;

  \fill[green!30] (axis cs:1,1) circle[radius=0.45];
\end{scope}

\addplot[
    blue,thick,
    domain=0:360,
    samples=300
]
({1+0.45*cos(x)},{1+0.45*sin(x)});

\addplot[
    red,thick,
    domain=0.4:1.167,
    samples=300
]
({x},{sqrt(max(0,2.7254-2*x^2))});

\node[circle,fill=orange,inner sep=2pt,
label=above right:{$x^0=(1.3,0.7)$}]
at (axis cs:1.2,0.7) {};

\node[circle,fill=black,inner sep=2pt,
label=below left:{$x^*$}]
at (axis cs:0.6818,0.6818) {};

\node[circle,fill=green!60!black,inner sep=2.5pt,
label=above:{$x^{5000}$}]
at (axis cs:0.6813,0.6823) {};

\node at (axis cs:0.9,0.9) {$Q$};
\node[blue] at (axis cs:1.25,1.38) {$\partial M$};
\node[red] at (axis cs:0.8,1.3) {$2x_1^2+x_2^2=2.7254$};

\end{axis}
\end{tikzpicture}
\caption{Feasible region $Q$, initial point $x^0$,
exact solution $x^*$, and final iterate $x^{5000}$.}
\label{fig:optimization_feasible_region}
\end{figure}
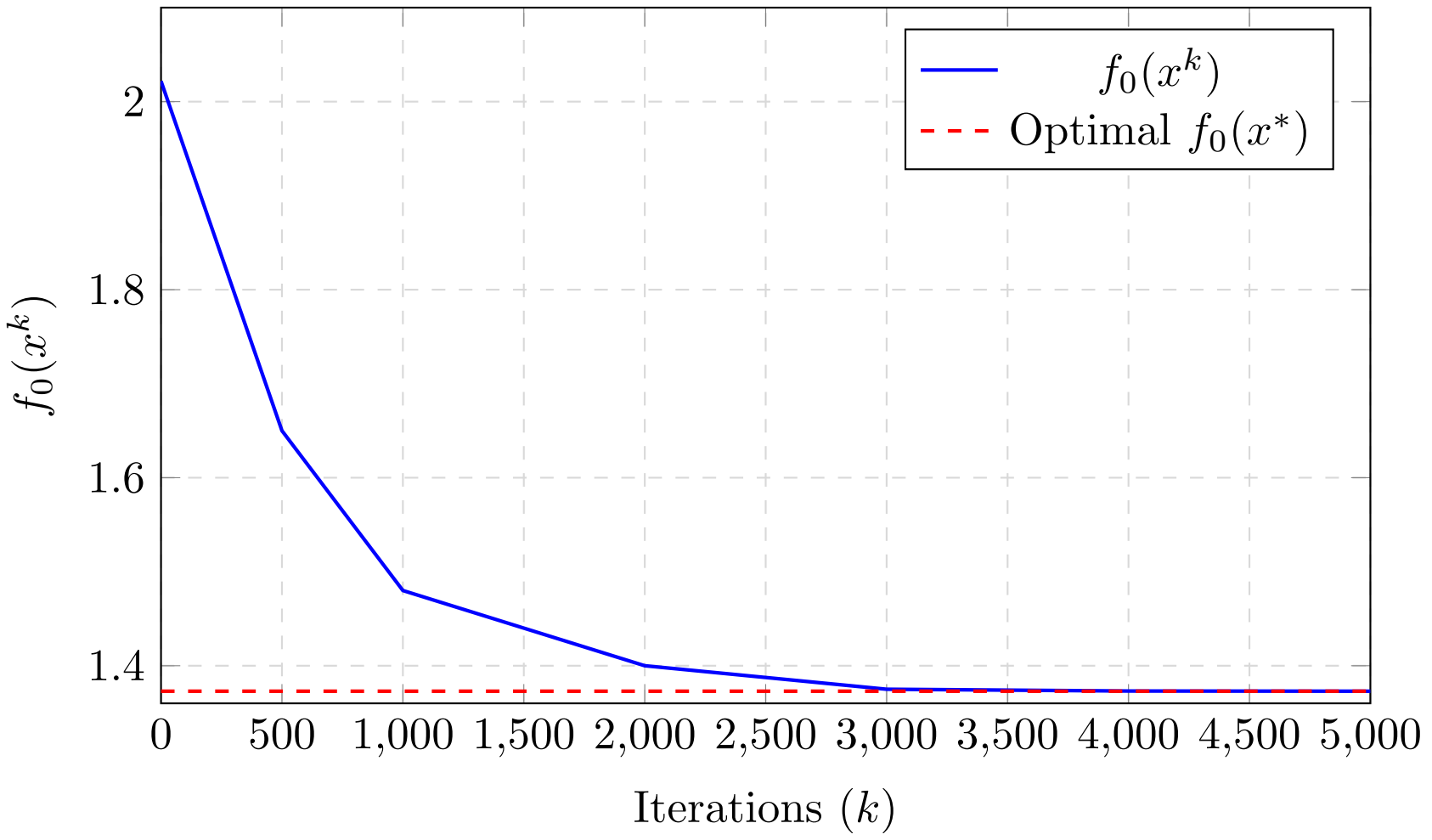
\begin{figure}[H]
\centering
\begin{tikzpicture}
\begin{axis}[
    xlabel = {Iterations ($k$)},
    ylabel = {$f_0(x^k)$},
    xmin = 0, xmax = 5000,
    ymin = 1.36, ymax = 2.1,
    grid = major,
    grid style = {dashed, gray!30},
    width = 11cm,
    height = 7.cm,
    legend pos = north east
    ]

    \addplot[blue, thick] coordinates {
        (0, 2.0218)
        (500, 1.65)
        (1000, 1.48)
        (2000, 1.40)
        (3000, 1.375)
        (4000, 1.3731)
        (5000, 1.37290178)
    };
    \addlegendentry{$f_0(x^k)$}

    \addplot[red, dashed, thick] coordinates {
        (0, 1.37290097)
        (5000, 1.37290097)
    };
    \addlegendentry{Optimal $f_0(x^*)$}

\end{axis}
\end{tikzpicture}
\caption{Convergence of the objective function $f_0(x^k)$ over 5000 iterations.}
\label{fig:objective_convergence}
\end{figure}

\begin{remark}[Practical performance and adaptive switching]
\label{rem:practical_switching}
In practical applications, estimating the problem-dependent constants $L_\Omega$ and $\gamma$ may be challenging, leading to conservative theoretical bounds for the parameter $\widetilde R$, such as $\widetilde R \geq 105$. We therefore examine the empirical behavior of Algorithm~\ref{alg:stochastic_switching} for a moderate value of $\widetilde R=1.0$.

For an initial constraint shift yielding $\widetilde f(x^0)\approx 2.4614$, the algorithm is run for $N=1000$ iterations. During iterations $k=0,\ldots,17$, the condition $\widetilde f(x^k)>\rho_k\widetilde R$ holds, and the algorithm selects the constraint-reduction direction
$p^k=g_{\xi_k}(x^k,\alpha_k,u_k).$
Starting at iteration $k=18$, the condition $\widetilde f(x^k)\leq\rho_k\widetilde R$ is satisfied, triggering a switch to the objective-minimization direction $ p^k=g_0(x^k,\alpha_k,u_k).$

After $N=1000$ iterations, the algorithm yields $x^{1000}\approx(0.6815,0.6821)$, with $f_0(x^{1000})\approx1.3730$ and $f_1(x^{1000})\approx-0.5370$. These results illustrate the adaptive switching mechanism and demonstrate that the algorithm attains a feasible point in this numerical experiment, despite using a value of $\widetilde R$ substantially below the theoretical bound.
\end{remark}

\section{Declarations}








\subsection{Author Contributions}

Both authors contributed equally to the study conception, design and implementation and wrote and corrected the manuscript.

\subsection{Competing Interests}

There are no conflicts of interest or competing interests related to this manuscript.

\subsection{Funding}

This research was partially supported by ANID--Chile under project Fondecyt Regular 1241040.

\subsection{Acknowledgements}
The authors wish to thanks the reviewers for their comments and remarks that helped to improve the quality of the paper.

\end{document}